\documentclass[runningheads,a4paper]{llncs}
\usepackage{a4wide}

\usepackage{amsmath,amsfonts,amssymb}
\usepackage{ifthen}
\usepackage{array}
\usepackage{url}
\usepackage{multirow}
\usepackage{graphicx}
\usepackage{enumitem}

\usepackage{graphicx}

\usepackage{url}
\usepackage{tabularx}
\usepackage{longtable}
\usepackage{pdflscape}

\RequirePackage[T1]{fontenc}
\RequirePackage[utf8]{inputenc}

\usepackage{soul} 
\usepackage[usenames,dvipsnames,svgnames,table]{xcolor}  
\usepackage{textcase}

\usepackage{multicol}

\newsavebox\ltmcbox

\newcounter{entryno}
\urldef{\mailsa}\path|oguz.yayla@ricam.oeaw.ac.at|

\newcommand{\keywords}[1]{\par\addvspace\baselineskip
\noindent\keywordname\enspace\ignorespaces#1}

\usepackage{setspace}

\newenvironment{pf}{\noindent\textbf{Proof.}}{\hfill{$\Box$}}
 
\usepackage{lipsum}

\newcommand\blfootnote[1]{%
  \begingroup
  \renewcommand\thefootnote{}\footnote{#1}%
  \addtocounter{footnote}{-1}%
  \endgroup
}

\begin{document}


\newcommand{\binomial}[2]{\left(\begin{array}{c}#1\\#2\end{array}\right)}
\newcommand{\zar}{{\rm zar}}
\newcommand{\an}{{\rm an}}
\newcommand{\red}{{\rm red}}
\newcommand{\codim}{{\rm codim}}
\newcommand{\rank}{{\rm rank}}
\newcommand{\Pic}{{\rm Pic}}
\newcommand{\Div}{{\rm Div}}
\newcommand{\Hom}{{\rm Hom}}
\newcommand{\im}{{\rm im}}
\newcommand{\Spec}{{\rm Spec}}
\newcommand{\sing}{{\rm sing}}
\newcommand{\reg}{{\rm reg}}
\newcommand{\Char}{{\rm char}}
\newcommand{\Tr}{{\rm Tr}}
\newcommand{\tr}{{\rm tr}}
\newcommand{\supp}{{\rm supp}}
\newcommand{\Gal}{{\rm Gal}}
\newcommand{\Min}{{\rm Min \ }}
\newcommand{\Max}{{\rm Max \ }}
\newcommand{\Span}{{\rm Span  }}

\newcommand{\Frob}{{\rm Frob}}
\newcommand{\lcm}{{\rm lcm}}

\newcommand{\soplus}[1]{\stackrel{#1}{\oplus}}
\newcommand{\dlog}{{\rm dlog}\,}
\newcommand{\limdir}[1]{{\displaystyle{\mathop{\rm
lim}_{\buildrel\longrightarrow\over{#1}}}}\,}
\newcommand{\liminv}[1]{{\displaystyle{\mathop{\rm
lim}_{\buildrel\longleftarrow\over{#1}}}}\,}
\newcommand{\boxtensor}{{\Box\kern-9.03pt\raise1.42pt\hbox{$\times$}}}
\newcommand{\sext}{\mbox{${\mathcal E}xt\,$}}
\newcommand{\shom}{\mbox{${\mathcal H}om\,$}}
\newcommand{\coker}{{\rm coker}\,}
\renewcommand{\iff}{\mbox{ $\Longleftrightarrow$ }}
\newcommand{\onto}{\mbox{$\,\>>>\hspace{-.5cm}\to\hspace{.15cm}$}}


\newcommand{\sA}{{\mathcal A}}
\newcommand{\sB}{{\mathcal B}}
\newcommand{\sC}{{\mathcal C}}
\newcommand{\sD}{{\mathcal D}}
\newcommand{\sE}{{\mathcal E}}
\newcommand{\sF}{{\mathcal F}}
\newcommand{\sG}{{\mathcal G}}
\newcommand{\sH}{{\mathcal H}}
\newcommand{\sI}{{\mathcal I}}
\newcommand{\sJ}{{\mathcal J}}
\newcommand{\sK}{{\mathcal K}}
\newcommand{\sL}{{\mathcal L}}
\newcommand{\sM}{{\mathcal M}}
\newcommand{\sN}{{\mathcal N}}
\newcommand{\sO}{{\mathcal O}}
\newcommand{\sP}{{\mathcal P}}
\newcommand{\sQ}{{\mathcal Q}}
\newcommand{\sR}{{\mathcal R}}
\newcommand{\sS}{{\mathcal S}}
\newcommand{\sT}{{\mathcal T}}
\newcommand{\sU}{{\mathcal U}}
\newcommand{\sV}{{\mathcal V}}
\newcommand{\sW}{{\mathcal W}}
\newcommand{\sX}{{\mathcal X}}
\newcommand{\sY}{{\mathcal Y}}
\newcommand{\sZ}{{\mathcal Z}}


\newcommand{\A}{{\mathbb A}}
\newcommand{\B}{{\mathbb B}}
\newcommand{\C}{{\mathbb C}}
\newcommand{\D}{{\mathbb D}}
\newcommand{\E}{{\mathbb E}}
\newcommand{\F}{{\mathbb{F}}}
\newcommand{\G}{{\mathbb G}}
\newcommand{\HH}{{\mathbb H}}
\newcommand{\I}{{\mathbb I}}
\newcommand{\J}{{\mathbb J}}
\newcommand{\M}{{\mathbb M}}
\newcommand{\N}{{\mathbb N}}
\renewcommand{\P}{{\mathbb P}}
\newcommand{\Q}{{\mathbb Q}}
\newcommand{\T}{{\mathbb T}}
\newcommand{\U}{{\mathbb U}}
\newcommand{\V}{{\mathbb V}}
\newcommand{\W}{{\mathbb W}}
\newcommand{\X}{{\mathbb X}}
\newcommand{\Y}{{\mathbb Y}}
\newcommand{\Z}{{\mathbb Z}}


\newcommand{\be}{\begin{eqnarray}}
\newcommand{\ee}{\end{eqnarray}}
\newcommand{\nn}{{\nonumber}}
\newcommand{\dd}{\displaystyle}
\newcommand{\ra}{\rightarrow}
\newcommand{\bigmid}[1][12]{\mathrel{\left| \rule{0pt}{#1pt}\right.}}
\newcommand{\cl}{${\rm \ell}$}
\newcommand{\clp}{${\rm \ell^\prime}$}


\mainmatter  

\title{Nearly perfect sequences with arbitrary out-of-phase autocorrelation}

\titlerunning{Non-existence of nearly perfect sequences}
%
%
\author{O\u guz Yayla}
\authorrunning{Yayla}

\institute{Johann Radon Institute for Computational and Applied Mathematics (RICAM), Austrian Academy of Sciences, Altenberger Strasse 69, A-4040 Linz, Austria\\
\mailsa}

%
%

\toctitle{Lecture Notes in Computer Science}
\tocauthor{Authors' Instructions}
\maketitle

\begin{abstract}
In this paper we study nearly perfect sequences (NPS) via their connection to direct product difference sets (DPDS). We prove the connection between a $p$-ary NPS of period $n$ and type $\gamma$ and a cyclic $(n,p,n,\frac{n-\gamma}{p}+\gamma,0,\frac{n-\gamma}{p})$-DPDS for an arbitrary integer $\gamma$. Next, we present the necessary conditions for the existence of a $p$-ary NPS of type $\gamma$. We apply this result for excluding the existence of some $p$-ary NPS of period $n$ and type $\gamma$ for $n \leq 100$ and $\vert \gamma \vert \leq 2$.  We also prove the similar results for an almost $p$-ary NPS of type $\gamma$. Finally, we show the non-existence of some almost $p$-ary perfect sequences by showing the non-existence of equivalent cyclic relative difference sets by using the notion of multipliers.
\end{abstract}

\section{Introduction}

Let $\underline{a} = (a_0,a_1,\ldots,a_{n-1},\ldots)$ be a sequence of period $n$ with entries $a_i \in \C$. For $0\leq t \leq n-1$, the \textit{autocorrelation 
function} $C_{\underline{a}}(t)$ is defined by
$$ C_{\underline{a}}(t) = \sum_{i=0}^{n-1}{a_i\overline{a_{i+t}}},$$
where $\overline{a}$ is the complex conjugate of $a$. The values $C_{\underline{a}}(t)$ at $1 \leq t \leq n-1$ are called \textit{the out-of-phase autocorrelation coefficients} of $\underline{a}$.

Let $p$ 
be a prime and $\zeta_p\in \C$ be a primitive $p$-th root of unity. If $a_i = \zeta_p ^{b_i}$ for some integer $b_i$, with $i =0, 1, \ldots,n-1$, 
then $\underline{a}$ is called a \textit{$p$-ary sequence}. If $a_0 = 0$ and all other entries are a power of $\zeta_p$, then $\underline{a}$ is called an 
\textit{almost $p$-ary sequence}.


An (almost) $p$-ary sequence $\underline{a}$ of period $n$ is called \textit{perfect sequence} (PS) if all out-of-phase autocorrelation coefficients of $\underline{a}$ are 0. 
Similarly, an (almost) $p$-ary sequence $\underline{a}$ of period $n$ is called a \textit{nearly perfect 
sequence} (NPS) of type $\gamma$ if all out-of-phase autocorrelation coefficients of $\underline{a}$ are a constant $\gamma$. We write a NPS of type $\gamma = 0$ to denote a PS. (We also note that there is another notion of \textit{almost perfect sequences} which is a $p$-ary sequence $\underline{a}$ of period $n$ having $C_{\underline{a}}(t)=0$ for all $1\leq t \leq n-1$ -with exactly one exception. In this paper, we will study $p$-ary NPS and almost $p$-ary NPS.)

There are some applications of (almost) $p$-ary sequences of period $n$ if $C_{\underline{a}}(t)$ is small for $1 \leq t \leq n-1$. 
We refer to  \cite{BJL} and the references therein for such applications.

Nearly perfect sequences have been studied widely by many authors. For instance, Jungnickel and Pott \cite{JuPo99} studied binary nearly perfect sequences of type $\vert \gamma \vert \leq 2$, and gave existence and non-existence cases.  Ma and Ng \cite{MaNg2009} studied $p$-ary nearly perfect sequences of type $\vert \gamma \vert \leq 1$, and determined their existence status by using their connection to direct product difference sets (DPDS). Later Chee et al.~\cite{CTZ2010} extended the methods due to Ma and Ng \cite{MaNg2009} 
to almost $p$-ary nearly perfect sequences of  types $\gamma = 0$ and $\gamma = -1$. Then, \"Ozbudak et al.~\cite{OYY2012} proved the non-existence of almost $p$-ary perfect sequences 
at certain values. 
Recently, Winterhof et al.~\cite{WYZ2014} studied the existence of (almost) $m$-ary NPS for an integer $m$ via their connection to Butson-Hadamard matrices, certain Diophantine equations and ideal decomposition.

In this paper we prove a general equality between an (almost) $p$-ary (nearly) perfect sequence of type $\gamma$ for an arbitrary $\gamma \in \Z$ and a direct product difference set. 
Firstly, we prove the connection between a $p$-ary NPS of type $\gamma$ and a
DPDS for an arbitrary integer $\gamma$ (see Theorem \ref{th:DPDS}). By this result we prove necessary conditions for the existence of a $p$-ary NPS of type $\gamma$ (see Theorem \ref{th:main1}). 
Then we demonstrate the pairs $(n,p)$ such that the existence of a $p$-ary NPS of period $n$ and type $\gamma$ is excluded by Theorem \ref{th:main1} for $n \leq 100$ and $ \gamma = -1,0,1,2$. 
In particular, we exclude the existence of a 23-ary NPS of period 45 and type $\gamma = -1$ by Theorem \ref{th:main1}, which was an undecided case in \cite{MaNg2009}. 
We note that an (almost) $p$-ary NPS of period $n$ and type $\gamma$ does not exists if $n \geq 3$ and $\gamma \leq -2$ (see Lemma \ref{lem:det} below).

Next, we prove the counterpart of the results on a $p$-ary NPS of type $\gamma$ for an almost $p$-ary NPS of type $\gamma$ where $\gamma \in \Z$ by considering its connection to a DPDS (see Theorems \ref{th:DPDS:almost} and \ref{th:main2}). And, we demonstrate the pairs $(n,p)$ such that the existence of an almost $p$-ary NPS of period $n+1$ and type $\gamma$ are excluded by Theorem \ref{th:main2} for $n \leq 100$ and $ \gamma = -1,0,1,2$. 
In particular, we exclude the existence of a 7-ary NPS of period 77 and type $\gamma = -1$ by Theorem \ref{th:main2}, which was an undecided case in \cite{CTZ2010}. Furthermore, we present a generalization of these results for showing the non-existence of an almost $p$-ary NPS with $s \geq 1$ zero-symbols. 

Finally, we show the non-existence of certain almost $p$-ary PS of period $n$ via showing the non-existence of a regular $(n+1,p,n,\frac{n-1}{p})$-RDS by using the notion of multipliers (see Theorem \ref{th:100-3}). 

The paper is organized as follows. In Section \ref{sec.pre}, we present the definition of DPDS and preliminary results that we use later. We present our result on $p$-ary NPS of type $\gamma$ in Section \ref{sec:main}, and on almost $p$-ary NPS of type $\gamma$ in Section \ref{sec:almost}. Finally, we present a result by using the notion of multipliers in Section \ref{sec.mult}.

\section{Preliminaries}
\label{sec.pre}
We begin with the definition of a direct product difference set \cite{MaNg2009}.

\begin{definition}
Let $G = H \times N$, where the order of $H$ and $N$ are $m$ and $n$. A subset $R$ of $G$ is called an $(m,n,k,\lambda_1,\lambda_2,\mu)$ direct product difference set (DPDS) in $G$ relative to $H$ and $N$ if both of the following statements hold:
\begin{itemize}
\item[(i.)] $|R| = k$,
\item[(ii.)] Differences $r_1r_2^{-1}$, $r_1, r_2 \in R$ with $r_1 \neq r_2$ represent
\begin{itemize}
\item all non identity elements of $H$ exactly $\lambda_1$ times,
\item all non identity elements of $N$ exactly $\lambda_2$ times,
\item all non identity elements of $G \backslash H \cup N$ exactly $\mu$ times.
\end{itemize}
\end{itemize}
\end{definition}

We can also define a difference set by using the group-ring algebra notation. 
Let $\sum_{g \in R}{g} \in \Z[G]$ be an element of the group ring $\Z[G]$, for simplicity we will denote the sum by $R$. If $R$ is an $(m,n,k,\lambda_1,\lambda_2,\mu)$-DPDS in $G$ relative to $H$ and $N$ then 
\be \label{eqn:DPDS} 
RR^{(-1)} = (k - \lambda_1 - \lambda_2 + \mu) + (\lambda_1 - \mu)H + (\lambda_2 - \mu)P + \mu G
\ee
holds in $\Z[G]$.

We now present two known results that we will use in subsequent sections. The following lemma is noticed by Turyn \cite{Tur1965}. Let $q$ be a prime and $u = q^rw$ where $\gcd(q,w) = 1$. We say that $q$ is \textit{self-conjugate} modulo $u$ if $q^j \equiv -1 \mod w$ for some integer $j$.

\begin{lemma} \label{lem:turyn}
If $q$ is self conjugate modulo $u$, then $\overline{Q} = Q$ for any prime ideal divisor $Q$ of $q\Z[\zeta_u]$.
\end{lemma}

Next result is known as Ma's Lemma \cite{MaNg2009}. We denote by $L^\perp$ the subset of the character group which is principal on $L$.

\begin{lemma} \label{lem:ma}
Let $q$ a be a prime and $\alpha$ be a positive integer. Let $K$ be an abelian group such that either $q$ does not divide $|K|$ or the Sylow $q$-subgroup of $K$ is cyclic. Let $L$ be any subgroup of $K$ and $Y \in \Z[K]$ where coefficients of $Y$ lie between $a$ and $b$ where $a < b$. Suppose that
\begin{enumerate}
\item $q$ is self conjugate modulo $\exp(K)$,
\item $q^r \mid \chi(Y)\overline{\chi(Y)}$ for all $\chi \not\in L^\perp$ and $q^{r+1} \nmid \chi(Y)\overline{\chi(Y)}$ for some $\chi \not\in L^\perp$,
\item $\chi(Y) \neq 0$ for some $\chi \not\in L^\perp \cup Q^\perp$ where $Q = K$ if $q \nmid |K|$, and $Q$ is the subgroup of $K$ of order $q$ otherwise.
\end{enumerate}
Then
\begin{enumerate}
\item if $q \nmid |K|$, $r$ is even and $q^{r/2} \leq b-a$, 
\item if Sylow $q$-subgroup of $K$ is cyclic, $q^{\lfloor \frac{r}{2} \rfloor} \leq 2(b-a)$ when $L$ is a proper subgroup of $|K|$ and $q^{\lfloor \frac{r}{2} \rfloor} \leq b-a$ when $L=K$.
\end{enumerate}
\end{lemma}

In the last part of this section we give a known result on the non-existence of NPS of type $\gamma$ for $\gamma \leq -2$. This is a direct consequence of \cite[Theorem 2.5]{Brock}, see also \cite[Corollary 3.1]{WYZ2014}. We give a short proof below.

\begin{lemma} \label{lem:det}
Let $p$ be a prime number, $n \in \Z^+$ and $\gamma \in \Z$ such that $n\geq 2$ and $\gamma \leq 2$. Then a $p$-ary (almost) NPS of period $n$ and type $\gamma$ does not exist except the existence of a binary NPS of period 2 and type -2.
\end{lemma}

\begin{pf}
Assume the existence of a $p$-ary NPS of period $n \geq 2$  and type $\gamma \leq -2$, say  $\underline a =(a_0, a_1, \ldots , a_{n-1})$. Let $H = (h_{i,j})$
be a circulant matrix, that is $h_{i+1,j+1} = h_{i,j}$ for all $i, j$, defined
by $h_{0,j} = a_j$ for $j=0,1,\ldots,n-1$ then $H$ is a circulant near Butson-Hadamard matrix of order
$n$ satisfying $H\overline{H}^T= (n-\gamma)I + \gamma J$, where $I$ is the identity matrix and $J$ is the all 1 square matrix of order $n$. Hence, $\det(H\overline{H}^T)=((\gamma+1)n-\gamma)(n-\gamma)^{n-1}$. Since $\det(H\overline{H}^T)$ is a non-negative number, we obtain a contradiction. We finally note that $(-1,1)$ is a binary NPS of period 2 and type -2. The proof of the non-existence of an almost $p$-ary NPS of period $n \geq 2$  and type $\gamma \leq -2$ is very similar.
\end{pf}

\section{$p$-ary (nearly) perfect sequences} \label{sec:main}

In the following result we prove the general connection between a $p$-ary sequence of type $\gamma$ and a DPDS for an integer $\gamma$. This is a generalization of \cite[Theorems 4.2 and 5.1]{MaNg2009}) to an arbitrary integer $\gamma$. 

\begin{theorem} \label{th:DPDS}
Let p be a prime, $n \geq 2$ be an integer, and $\underline{a} = (a_0,a_1,\ldots,a_{n-1},\ldots)$ be a $p$-ary sequence of period $n$.
Let $H = \langle h \rangle$ and $P = \langle g \rangle $ be the (multiplicatively written) cyclic groups of order $n$ and $p$, respectively. Let $G$ be the group defined as $G=H \times P$.
We choose a primitive $p$-th root of 1, $\zeta_p \in \C$. For $0 \leq i \leq n-1$ let $b_i$ be the integer in $\{0,1,2,\ldots,p-1\}$ such that $a_i = \zeta_p^{b_i}$.
Let $R$ be the subset of $G$ defined as
\be \nn
R = \{(g^{b_i}h^i) \in G : 0 \leq i \leq n-1\}.
\ee
Then 
$\underline{a}$ is a $p$-ary NPS of type $\gamma$ if and only if $R$ is an $(n,p,n,\frac{n-\gamma}{p}+\gamma,0,\frac{n-\gamma}{p})$-DPDS in $G$ relative to $H$ and $N$. In particular, $p$ divides $n-\gamma$.

\end{theorem}
\begin{pf}
Let $A = \sum_{i = 0}^{n-1}{a_ih^i} \in \C[H]$. Then we have 
$$A\overline{A}^{(-1)} = \sum_{t = 0}^{n-1}{C_a(t)h^t}.$$
Let $\chi$ be a character on $P$. We extend $\chi$ to $G$ such that $\chi(h) = h$. Let $\sigma \in \Gal(\Q(\zeta_p)\backslash \Q)$ such that $\sigma(\zeta_p) = \chi(\zeta_p)$. 
If $\chi$ is a nonprincipal character on $P$, then we have $\chi(R) = A^{\sigma}$, and so
$$\chi(RR^{(-1)}) = (A\overline{A}^{(-1)})^{\sigma}.$$
On the other hand, if $\chi$ is a principal character on $P$, then we have 
$$\chi(R) = H.$$
Then
\be \nn
\chi(RR^{(-1)}) = \left\lbrace \begin{array}{ll} nH & \mbox{if } \chi \mbox{ is principal on }P\\ \sum_{t = 0}^{n-1}{C_a(t)^\sigma h^t} & \mbox{if } \chi \mbox{ is nonprincipal on $P$ } \end{array} \right.
\ee
If $a$ is a NPS of type $\gamma$, then 
\be \nn
\chi(RR^{(-1)}) = \left\lbrace \begin{array}{ll} nH & \mbox{if } \chi \mbox{ is principal on }P\\ n - \gamma + \gamma H & \mbox{if } \chi \mbox{ is nonprincipal on }P \end{array} \right.
\ee
By extending $\chi$ to $H$ we obtain
\be \nn
\chi(RR^{(-1)}) = \left\lbrace 
\begin{array}{ll} 
n^2 & \mbox{if } \chi \mbox{ is principal on }P  \mbox{ and } H\\ 
0 & \mbox{if } \chi \mbox{ is principal on }P  \mbox{ and nonprincipal on }H\\
(\gamma + 1)n -\gamma & \mbox{if } \chi \mbox{ is nonprincipal on }P   \mbox{ and principal on }H\\
n -\gamma & \mbox{if } \chi \mbox{ is nonprincipal on }P  \mbox{ and nonprincipal on }H\\ 
\end{array} 
\right.
\ee
On the other hand it is easy to see by using (\ref{eqn:DPDS}) that the same diagram holds for an $(n,p,n,\frac{n-\gamma}{p}+\gamma,0,\frac{n-\gamma}{p})$-DPDS for any character $\chi$ on $G$. Therefore we are done.
\end{pf}

Theorem \ref{th:DPDS} gives a way of showing the existence and the non-existence of NPS via using DPDS. We will use this method in proving the following theorem. For integers $q,r$ and $n$ we use $q^r||n$ to denote that $q^r|n$ but $q^{r+1} \nmid n$.

\begin{theorem} \label{th:main1}
Let $p$ be prime number, $n \in \Z^+$ and $\gamma \in \Z$ such that $|\gamma| < n$ and $p\mid n - \gamma$.  Suppose that there exists a $p$-ary NPS of type $\gamma$ and period $n$.
\begin{itemize}
\item[(i)]  For $\gamma = 0$, let $q \neq p$ be prime number dividing $n$ such that $q^r||n$ and $q$ is self-conjugate modulo $up$ for some divisor $u \geq 1$ of $n$. Then $r$ is even and $\dd q^{r/2} \leq \frac{n}{u}$.
\item[(ii)]  For $\gamma \neq 0$, let $q \neq p$ be prime number dividing $n- \gamma$ such that $q^r||n - \gamma$ and $q$ is self-conjugate modulo $up$ for some divisor $u>1$ of $n$. If $q \nmid u$, then $r$ is even and $\dd q^{r/2} \leq \frac{n}{u}$. If $q \mid u$ then $\dd q^{\lfloor r/2 \rfloor} \leq 2\frac{n}{u}$. 
\item[(iii)] For $\gamma \neq 0$, if $p^r||n - \gamma$ and $p$ is self-conjugate modulo $u$ for some divisor $u>1$ of $n$ such that $p \nmid u$, then $\dd p^{r/2} \leq 2\frac{n}{u}$ in case $r$ is even, and $\dd p^{(r+1)/2} \leq 4\frac{n}{u}$ in case $r$ is odd. 
\item[(iv)]  For $\gamma \neq 0$, let $q \neq p$ be prime number dividing $(\gamma +1)n - \gamma$ such that $q^r||(\gamma + 1)n - \gamma$ and $q$ is self-conjugate modulo $up$ for some divisor $u \geq 1$ of $n$. If $q \nmid u$, then $r$ is even. 
\end{itemize}
\end{theorem}
\begin{pf}
We note that (i) is already proved in \cite[Theorem 4.11]{MaNg2009}. Therefore, it is enough to prove (ii), (iii) and (iv). By using Theorem \ref{th:DPDS}, the existence of a p-ary NPS of period $n$ and  type $\gamma$ implies that there is an $(n,p,k,\frac{n-\gamma}{p}+\gamma,0,\frac{n-\gamma}{p})$-DPDS in $G = H \times P$ relative to $H = \langle h \rangle$ and $P = \langle g \rangle$ where $o(h) = n$ and $o(g) = p$.

Let $\rho : G \rightarrow K:= G / \langle h^u \rangle$ be the natural epimorphism. By using (\ref{eqn:DPDS}) with $(n,p,k,\frac{n-\gamma}{p}+\gamma,0,\frac{n-\gamma}{p})$-DPDS, we have
\be \nn
\rho(RR^{(-1)}) = (n - \gamma) + \gamma \frac{n}{u} \rho(H) - (\frac{n - \gamma}{p}) \rho(P) + (\frac{n - \gamma}{p}) \frac{n}{u} K
\ee
The coefficients of $\rho(R)$ lie between  0 and $\dd\frac{n}{u}$. 
Let $\chi$ be a nonprincipal character of $K$. Then
\be \label{diagram}
\chi(\rho(RR^{(-1)})) = \left\lbrace 
\begin{array}{ll} 
0 & \mbox{if } \chi \mbox{ is principal on }\rho(P)  \mbox{ and nonprincipal on }\rho(H)\\
(\gamma + 1)n -\gamma & \mbox{if } \chi \mbox{ is nonprincipal on }\rho(P)   \mbox{ and principal on }\rho(H)\\
n -\gamma & \mbox{if } \chi \mbox{ is nonprincipal on }\rho(P)  \mbox{ and nonprincipal on }\rho(H)\\ 
\end{array} 
\right.
\ee
Assume that $q^r||n-\gamma$ and $q$ is self-conjugate modulo $\exp(K) = up$ for some divisor $u>1$  of $n$. We will use Lemma \ref{lem:ma} with $L = \rho(H)$ and $Y = \rho(R)$ to prove (ii). We note that  $\rho(H)$ is nontrivial as $u >1$.  Then we have $q^r \mid \chi(Y)\overline{\chi(Y)}$ for all $\chi \not\in L^\perp$ and $q^{r+1} \nmid \chi(Y)\overline{\chi(Y)}$ for $\chi \not\in L^\perp \cup \rho(P)^\perp$. If $q \mid |K|$, then it is clear that $\chi(Y) \neq 0$ for some $\chi \not\in \rho(Q)^\perp \cup \rho(P)^\perp$, where $Q$ is the Sylow $q$-subgroup of $K$ since $\chi(Y)\overline{\chi(Y)} = n - \gamma$. Then by Lemma \ref{lem:ma} we get (ii). 

Similarly, if $q = p$ and $p \nmid u$ then the Sylow $p$-subgroup of $K$ is cyclic. When $r$ is even, we obtain the first part of (iii) by using Lemma \ref{lem:ma}.  When $r$ is odd, we define 
$$R' := R\sum_{t = 1}^{p-1}{\left( \frac{t}{p} \right)g^t}$$
where $\left(\dd\frac{t}{p} \right)$ denotes the Legendre symbol. We note that the coefficients of $R'$ are $-1,0,+1$, so the coefficients of $\rho(R')$ are between $-n/u$ and $n/u$. We also note that $$\sum_{t = 1}^{p-1}{\left( \frac{t}{p} \right)g^t} \ \overline{\sum_{t = 1}^{p-1}{\left( \frac{t}{p} \right)g^t}} = p$$ by using \cite[Lemma 4.5]{MaPo1995}. And so, $p^{r+1} \mid \chi(\rho(R'))\overline{\chi(\rho(R'))}$ for all $\chi \not\in L^\perp$ and $p^{r+2} \nmid \chi(\rho(R'))\overline{\chi(\rho(R'))}$ for $\chi \not\in L^\perp \cup \rho(P)^\perp$. Thus by using Lemma \ref{lem:ma} we prove the second part of (iii).

Finally, it is clear from (\ref{diagram}) that there exists a nonprincipal character $\chi$ of $K$ such that $\chi(Y)\overline{\chi(Y)} =(\gamma + 1)n - \gamma$. If there exists a prime number $q \neq p$ such that $q^r||(\gamma + 1)n - \gamma$ and $q$ is self-conjugate modulo $up$ for some divisor $u \geq 1$ of $n$, then $r$ must be even by Lemma \ref{lem:turyn}. This proves (iv).
\end{pf}

Ma and Ng \cite{MaNg2009} present tables showing the existence status of $p$-ary (nearly) perfect sequences of period $n$ for $|\gamma| \leq 1$ and $2 \leq n \leq 50$.
We extend the tables given in \cite{MaNg2009} to $|\gamma| \leq 2$ and $2 \leq n \leq 100$. In addition, we update the tables for some undecided cases. Below, we explain existence, non-existence and undecided cases. 
We present the detailed tables in Appendix. The empty rows in the tables are undecided cases. 
The case $p=2$ is extensively studied in \cite{JuPo99}, therefore in this section we only deal with the case that $p$ is an odd prime.

For $\gamma = 0$, it is known that an $(n,p,n,n/p,0,n/p)$-DPDS in $\Z_n \times \Z_p$ relative to $\Z_n$ and $\Z_p$ exists for $n = p$ and $n = p^2$ where $p$ is an odd prime (see \cite[Theorem 2.2.9]{pott1995} and \cite[Theorem 2.3]{MaSc1995} respectively). Therefore, a $p$-ary PS of period n for $n=p$ and $n=p^2$ exist. For $n \leq 100$, Theorem \ref{th:main1} excludes the existence at all other pairs $(n,p)$ except a few undecided cases $(n,p) \in \{(28,7),(33,11),(39,13),(55,11),(56,7),(63,3),(69,23),(84,3),(92,23),(95,19),(99,11) \}$.

For $\gamma = -1$, it is known that an $(n,p,n,(n+1-p)/p,0,(n+1)/p)$-DPDS in $\Z_n \times \Z_p$ relative to $\Z_n$ and  $\Z_p$ exists for $n = q-1$ where $q$ is a power of $p$ (see \cite{HK1998}). Therefore, a $p$-ary NPS of period $q-1$ for  $q$ is a power of $p$ exists. For $n \leq 100$, Theorem \ref{th:main1} excludes the existence at all other pairs $(n,p)$ except $(n,p) \in  \{$(19,5), (23,3), (27,7), (32,11), (35,3), (38,3), (41,3), (44,5), (47,3), (55,7), (56,3), (59,5), (65,3), (65,11), (67,17), (71,3),  (73,37), (74,3), (76,7), (79,5), (83,7), (91,23), (92,3), (93,47), (98,11), (99,5)$\}$. We note that 23-ary NPS of period 45 for $\gamma = -1$ was an undecided case in \cite{MaNg2009}. Theorem \ref{th:main1} (iii) with $q = 23$ and $u = 9$ shows that such a sequence does not exist. It is assumed in \cite{MaNg2009} that Theorem \ref{th:main1} (ii) holds for $\gamma = -1$ and $u = 1$. However, we show in the proof that it does not hold for $u =1$. Because of this reason, we say that the cases (23,3) and (41,3) are undecided which were given to be nonexistent in \cite{MaNg2009}.

For $\gamma = 1$, Theorem \ref{th:main1} excludes the existence of a $p$-ary NPS of period $n$ for the pairs $(n,p)$ such that $n \leq 100$ except the ones $(n,p) \in  \{$(5,3), (13,3), (15,7), (22,3), (25,3), (27,13), (31,3), (31,5), (40,3), (40,13), (45,11), (49,3), (51,5), (56,11), (57,7), (63,31), (64,3), (64,7), (70,23), (76,3), (79,3), (79,13), (85,3), (85,7), (95,47), (96,19), (97,3), (99,7), (100,11)$\}$. And, it is known that a $p$-ary NPS for $\gamma = 1$ of period $n$ exists for pairs $(n,p) \in  \{$(5,3), (13,3)$\}$ (see \cite{MaNg2009}).

For $\gamma = 2$, Theorem \ref{th:main1} excludes the existence of a $p$-ary NPS of period $n$  for the pairs $(n,p)$ such that $n \leq 100$ except the ones $(n,p) \in \{$(5,3), (9,7), (11,3), (16,7), (17,3), (22,5), (23,3), (26,3), (30,7), (33,31), (35,3), (35,11), (37,7), (41,3), (50,3), (57,5), (58,7), (59,3), (59,19), (65,3), (77,3), (81,79), (86,3), (95,3), (98,3)$\}$. On the other hand, an exhaustive search says that no $p$-ary NPS of period $n$ and type $\gamma = 2$ exists for the pairs $(n,p) = (9,7)$, but it exists for the pairs $(n,p) \in  \{(5,3), (17,3)\}$: $(\zeta_3^2, \zeta_3^2, \zeta_3^2, \zeta_3^2, 1)$ is a 3-ary NPS of period 5 and type 2, and $(\zeta_3^2, \zeta_3^2, \zeta_3^2, 1,\zeta_3^2, 1, 1, \zeta_3, 1, 1,\zeta_3^2, 1, \zeta_3^2, \zeta_3^2,\zeta_3^2,1, 1)$ is a 3-ary NPS of period 17 and type 2.

\section{Almost $p$-ary (nearly) perfect sequences}
\label{sec:almost}

In the following we prove the equality between an almost $p$-ary sequence of type $\gamma$ and a DPDS for an integer $\gamma$. This is a generalization of \cite[Theorems 1 and 6]{CTZ2010} to an arbitrary integer $\gamma$. Its proof is similar to the proof of Theorem \ref{th:DPDS}.

\begin{theorem} \label{th:DPDS:almost}
Let p be a prime, $n \geq 2$ be an integer, and $\underline{a} = (a_0,a_1,\ldots,a_{n},\ldots)$ be an almost $p$-ary sequence of period $n+1$.
Let $H = \langle h \rangle$ and $P = \langle g \rangle $ be the (multiplicatively written) cyclic groups of order $n+1$ and $p$. Let $G$ be the group defined as $G=H \times P$.
We choose a primitive $p$-th root of 1, $\zeta_p \in \C$. For $1 \leq i \leq n$ let $b_i$ be the integer in $\{0,1,2,\ldots,p-1\}$ such that $a_i = \zeta_p^{b_i}$. We define $a_0 := 0$.
Let $R$ be the subset of $G$ defined as
\be \nn
R = \{(g^{b_i}h^i) \in G : 1 \leq i \leq n\}.
\ee
Then 
$\underline{a}$ is an almost $p$-ary NPS of type $\gamma$ if and only if $R$ is an $(n+1,p,n,\frac{n-\gamma-1}{p}+\gamma,0,\frac{n-\gamma-1}{p})$-DPDS in $G$ relative to $H$ and $P$. In particular, $p$ divides $n-\gamma-1$.
\end{theorem}
\begin{pf}
We define $$A:= \sum_{i = 0}^{n}{a_ih^i} \in \C[H],$$ where $a_0 =0$. By applying steps in the proof of Theorem \ref{th:main1} we complete the proof.
\end{pf}

Let $\chi$ be a character of $G = H \times P$ where $G$ is defined as in Theorem \ref{th:DPDS:almost}. If $R \in \Z[G]$ is an $(n+1,p,n,\frac{n-\gamma-1}{p}+\gamma,0,\frac{n-\gamma-1}{p})$-DPDS in $G$ relative to $H$ and $P$, then we have
\be \label{diag:almost}
\chi(RR^{(-1)}) = \left\lbrace 
\begin{array}{ll} 
n^2 & \mbox{if } \chi \mbox{ is principal on }P  \mbox{ and } H\\ 
1 & \mbox{if } \chi \mbox{ is principal on }P  \mbox{ and nonprincipal on }H\\
(\gamma + 1)n & \mbox{if } \chi \mbox{ is nonprincipal on }P   \mbox{ and principal on }H\\
n -\gamma & \mbox{if } \chi \mbox{ is nonprincipal on }P  \mbox{ and nonprincipal on }H\\ 
\end{array} 
\right.
\ee

We note that Lemma \ref{lem:ma} can  be applied to almost $p$-ary NPS only for $\gamma = 0$ and $\gamma = -1$. Because in other cases we can not find a prime dividing $\chi(RR^{(-1)})$ for any nonprincipal character defined over a subset of $G$, see (\ref{diag:almost}). In the cases $\gamma = 0$ or $\gamma = -1$ one can use Lemma \ref{lem:ma} with $L = P$, and these cases are already considered in \cite[Theorems 3 and 8]{CTZ2010}. We state them in Theorem \ref{th:main2} (i) and (ii), respectively. On the other hand, we extend Theorem \ref{th:DPDS} (iv) to almost $p$-ary NPS in Theorem \ref{th:main2} (iii).

\begin{theorem} \label{th:main2}
Let $p$ be prime number and $n \in \Z^+$, $\gamma \in \Z$ such that $|\gamma| < n$.  Suppose that there exists a type $\gamma$ almost $p$-ary NPS of period $n+1$.
\begin{itemize}
\item[(i)]  For $\gamma = 0$, let $q \neq p$ be prime number dividing $n$ such that $q^r||n$ and $q$ is self-conjugate modulo $up$ for some divisor $u \geq 1$ of $n+1$. Then $r$ is even and $\dd q^{r/2} \leq \frac{n+1}{u}$.
\item[(ii)]  For $\gamma = -1$, let $q \neq p$ be prime number dividing $n- \gamma$ such that $q^r||n - \gamma$ and $q$ is self-conjugate modulo $up$ for some divisor $u>1$ of $n+1$. If $q \nmid u$, then $r$ is even and $\dd q^{r/2} \leq \frac{n+1}{u}$. If $q \mid u$ then $\dd q^{\lfloor r/2 \rfloor} \leq 2\frac{n+1}{u}$. 
\item[(iii)] For $\gamma \neq 0$, let $q \neq p$ be prime number dividing $(\gamma +1)n$ $($or $n - \gamma)$ such that $q^r||(\gamma + 1)n$ $($or resp.~$q^r || n - \gamma)$ and $q$ is self-conjugate modulo $up$ for some divisor $u \geq 1$ $($or resp.~$u > 1)$ of $n+1$. If $q \nmid u$, then $r$ is even. 
\end{itemize}
\end{theorem}
\begin{pf}
The conclusions (i) and (ii) are already proved in \cite{CTZ2010}. For the proof of (iii) we use (\ref{diag:almost}). We have a nonprincipal character $\chi$ of $G$ such that $\chi(R)\overline{\chi(R)} =(\gamma + 1)n$. Hence if there exists a prime $q \neq p$ such that $q^r||(\gamma + 1)n$ and self-conjugate modulo $up$ for some divisor $u \geq 1$ of $n+1$, then by Lemma \ref{lem:turyn} $r$ must be even. Similarly, if there exists a prime $q \neq p$ such that $q^r||n - \gamma$ and self-conjugate modulo $up$ for some divisor $u > 1$ of $n+1$, then $r$ must be even.
\end{pf}

We note that Theorem \ref{th:DPDS:almost} (iii) can be extended to the sequences consisting of\linebreak $\underline{a} = (a_0,a_1,\ldots,a_{n+s-1},\ldots)$ of period $n+s$ with $a_{i_j} = 0$ for all $j=1,2,\ldots , s$ where $\{ i_1,i_2, \ldots , i_j \} \subset \{0,1, \ldots , n+s-1\}$ and $a_i = \zeta_p^{b_i}$ for some integer $b_i$, $i \in  \{0,1, \ldots , n+s-1\} \backslash \{ i_1,i_2, \ldots , i_j \}$ where $\zeta_p$ is a $p$-th root of unity  in $\C$. We call $\underline{a}$ an \textit{almost $p$-ary sequence with $s$ zero-symbols}. 
\begin{corollary}
Let $p$ be prime number and $n \in \Z^+$, $\gamma \in \Z$ such that $0<|\gamma| < n$.   Let $q \neq p$ be prime number dividing $(\gamma +1)n + (s-1)\gamma$ $($or $n - \gamma)$ such that $q^r||(\gamma + 1)n+ (s-1)\gamma$ $($or resp.~$q^r || n - \gamma)$ and $q$ is self-conjugate modulo $up$ for some divisor $u \geq 1$ $($or resp.~$u > 1)$ of $n+s$ and $q \nmid u$. If  there exists a type $\gamma$ almost $p$-ary NPS of period $n+s$ with $s$ zero-symbols, then $r$ is even. 
\end{corollary}

Chee et al.~\cite{CTZ2010} extend the results in \cite{MaNg2009} and present tables showing the existence status of almost $p$-ary NPS of period $n+1$ for $ \gamma = 0$ and $\gamma = 1$ and $2 \leq n \leq 100$.
We extend the tables given in \cite{CTZ2010} for $|\gamma| \leq 2$ and $2 \leq n \leq 100$. In addition, we update the tables in \cite{CTZ2010} for some undecided cases. Below, we explain existence, non-existence and undecided cases. 
We present the detailed tables in Appendix. The empty rows in the tables are undecided cases.

For $\gamma = 0$, it is known that an $(n+1,p,n,(n-1)/p,0,(n-1)/p)$-DPDS in $\Z_{n+1} \times \Z_p$ relative to $\Z_{n+1}$ and $\Z_p$ exists if $n$ is a prime power and $p$ divides $n-1$ (see \cite[Theorem 2.2.12]{pott1995}). Therefore, an almost $p$-ary PS of period $n+1$ exists when $n$ is a prime power and $p$ is a prime divisor of $n-1$. By using Theorem \ref{th:main2}  and results in \cite{CTZ2010,OYY2012} we obtain for $n \leq 100$ that an almost $p$-ary PS of period $n+1$ at all other cases do not exist except the undecided pairs  $(n,p) \in \{(63,31),(77,19),(91,3),(92,7),(93,23) \}$. In Section \ref{sec.mult}, we exclude the existence at the cases $(n,p) \in \{(63,31),(91,3),(92,7),(93,23) \}$ by using multipliers.

For $\gamma = -1$, it is known that an $(n+1,p,n,n/p,0,n/p)$-DPDS in $\Z_{n+1} \times \Z_p$ relative to  $\Z_{n+1}$ and $\Z_p$ exists for $n = q-1$ where $q$ is a prime and $p$ divides $q-1$ (see \cite[example 5.3.2]{pott1995}). Therefore, an almost $p$-ary NPS of period $q$ for  $q$ is a prime and $p \mid q-1$ exists. Theorem \ref{th:main1} excludes the existence of an almost $p$-ary NPS of period $n+1$ for the remaining pairs $(n,p)$ except $(n,p) \in \{$(20,5), (21,7), (26,13), (27,3), (34,17), (35,5), (38,2), (38,19), (44,11), (48,3), (50,5), (51,3), (54,2), (54,3), (63,7), (68,17), (75,3), (76,19), (84,3), (84,7), (90,3), (91,7),  (92,23), (93,31), (98,7), (99,11)$\}$. We note that 7-ary NPS of period 77 for $\gamma = -1$ was an undecided case in \cite{CTZ2010}. Theorem \ref{th:main2} (ii) for $q = 3$ and $u = 2$ shows that such a sequence does not exist. It is assumed in \cite{CTZ2010} that Theorem \ref{th:main1} (ii) holds for $\gamma = -1$ and $u = 1$. However, we show in the proof that it does not hold for $u =1$. Because of this reason, we say that the cases $(n,p) \in \{$(38,2), (38,19), (50,5), (54,2), (54,3), (68,17), (84,3), (84,7)$\}$  were undecided which are given to be nonexistent in \cite{CTZ2010}.

For $\gamma = 1$, Theorem \ref{th:main1} excludes the existence of an almost $p$-ary NPS of period $n+1$ for the pairs $(n,p)$ such that $n \leq 100$ except the cases $(n,p) \in \{$(8,2), (8,3), (9,7), (14,3), (16,7), (18,2), (22,5), (23,7), (24,11), (25,23), (26,3), (32,2), (32,3), (32,5), (37,7), (38,3), (40,19), (44,7), (46,11), (48,23), (49,47), (50,2), (50,3), (54,13), (58,7), (62,3), (62,5), (64,31), (72,5), (72,7), (73,71), (74,3), (81,79), (82,5), (88,43), (90,11), (94,23), (95,31), (96,47), (98,2), (98,3), (100,7)$\}$. In addition, we performed an exhaustive search for the pairs $(n,p) \in\{$(8,2), (8,3), (9,7), (14,3), (18,2), (32,2)$\}$, and we obtained that an almost $p$-ary NPS of period $n+1$ for $\gamma = 1$ exists for none of them.

For $\gamma = 2$, Theorem \ref{th:main1} excludes the existence of an almost $p$-ary NPS of period $n+1$ for the pairs $(n,p)$ such that $n \leq 100$ except the cases $(n,p)\in\{$(9,3), (12,3), (16,13),  (24,7), (25,11), (26,23), (27,2), (27,3), (29,13), (33,5), (36,3), (36,11), (39,3), (47,11), (48,3), (48,5), (49,23), (50,47), (60,19), (63,3), (66,7), (69,11), (72,23), (74,71), (75,2), (75,3), (81,3), (81,13), (84,3), (93,3), (96,31)$\}$. In addition, we exclude the existence of the pairs in $\{ ( 9, 3 ), (21,3), ( 27, 2 )\}$ by an exhaustive search. On the other hand, we have an example of almost $3$-ary NPS of period $13$ for  $\gamma = 2$: $(0,\zeta_3^2, \zeta_3^2, \zeta_3^2,  1,\zeta_3^2, \zeta_3, \zeta_3,\zeta_3^2, 1, \zeta_3^2, \zeta_3^2,\zeta_3^2)$.
   
\section{Non-existence by using multipliers}
\label{sec.mult}
An important method for the existence and the non-existence of some difference sets uses the notion of multiplier. In this section we prove the non-existence of almost $p$-ary PS at some values by showing the non-existence of the corresponding DPDS such that the existence of these values are not excluded by Theorem \ref{th:main2}.
We note that an $(n+1,p,n,\frac{n-\gamma-1}{p}+\gamma,0,\frac{n-\gamma-1}{p})$-DPDS in $G=H \times P$ relative to $H$ and $P$ for $\gamma =0$ is called an $(n+1,p,n,\frac{n-1}{p})$ \textit{relative difference set} (RDS) in $G=H \times P$ relative to $P$. 

Let $R$ be a subset in $G$. For an integer $t$, let $R^{(t)}$ denote the subset $R^{(t)}=\{r^t:r\in R\}\subset G$.
Assume that $gcd(t,|G|)=1$. $t$ is called a \textit{multiplier} of $R$ if there exists $g \in G$ such that
\be \nn
R^{(t)}=Rg=\{rg:r\in R\}\subset G.
\ee

There is a nice method for the existence and the non-existence of certain RDS that we recall here (see \cite{CTZ2010} page 406). Assume that $R$ is an $(m,n,k,\lambda)$-RDS in $G$ relative to $P$ such that $k^2\neq \lambda mn$ and $t$ is a multiplier of $R$. Let $\Omega$ be the set of orbits of $G$ under the action
$x \ra x^t$. Then, there exists a collection $\Phi$ of orbits (i.e. a subset $\Phi \subseteq \Omega$) such that
\be \nn
R=\bigsqcup_{A \in \Phi}{A},
\ee
where $A$ is an orbit in $\Phi$. This gives a strict condition on the existence and the non-existence of RDS.

In addition, we use the following result in proving the non-existence of RDS at certain values, see \cite[Lemma 5.4]{BJL} or \cite[Proposition 1]{OYY2012}.
\begin{lemma} \label{lem.core}
Let $R$ be an $(n+1,p,n,\frac{n-1}{p})$-RDS in $G=\Z_{n+1}\times \Z_p$ relative to $P=\Z_p$.
Let $R$ have $s_i$ many elements having $i$ in the second component for $i=0,1,2,\ldots,p-1$.
Then $$\sum_{j=0}^{p-1}{s_j^2}=\frac{n(n+p-1)}{p} \mbox{ and }  \sum_{j=0}^{p-1}{s_js_{j-i}}=\frac{n(n-1)}{p}$$
for each $i = 1,2,\ldots,\left\lceil(p-1)/2\right\rceil$, where subscripts are computed modulo $p$.
\end{lemma}

By using the method presented above and Lemma \ref{lem.core} we prove the following result.
\begin{theorem} \label{th:100-3}
There does not exist an almost $p$-ary perfect sequence with period $n+1$ for the pairs $(n,p) \in \{ (63,31), (91,3), (92,7), (93,23) \}$.
\end{theorem}
\begin{pf}
We present here the proof of the case (91,3). The others can be proven similarly. Assume that there exists an almost 3-ary PS of period 92. Using Theorem \ref{th:DPDS:almost} we have
an $(92,3,91,30)$-RDS $R$ in $\Z_{92} \times \Z_3$ relative to $\Z_3$. It is easy to see that $t=13$ is a multiplier of $R$.
Indeed let $\zeta$ be a primitive $276$-th root of $1$ in $\C$. We have $91=7\cdot 13$ and $\zeta^{13}=(\zeta^7)^{199}$. We tabulate the set $\Omega$ of orbits under the action $x \ra x^{13}$ in $G$ in Table \ref{tab.91-3}, and see that there are 12 orbits of length 1 and 24 orbits of length
11. Moreover, we may assume that there exists a subset $\Phi \subset \Omega$ satisfying
\be \nn
R=\bigsqcup_{A \in \Phi}{A}.
\ee
As in Lemma \ref{lem.core}, let $s_0,s_1$ and $s_2$ denote the number of elements in $R$ with the second component 0,1 and 2 respectively.
Using Lemma \ref{lem.core} we obtain that
\be \label{eqn.100}
s_0^2+s_1^2+s_2^2=\frac{91\cdot 93}{3} =2821.
\ee
We can choose orbits of length 1 in $\Omega$ with only the same second component. 
Thus we can only choose at most 4 orbits of length 1 for $\Phi$ covering $R$. Moreover, we can not choose two orbits $A_1$ and $A_2$ of length 11 of $\Omega$ such that with the same first components. Otherwise, difference of elements in $A_1$ and $A_2$ with the same first component gives an element in the forbidden group. We have 8 distinct subsets of orbits in $\Omega$ of length 11 with the same first component.

As $|R|=91$, it is clear from the lengths and the numbers of the orbits in Table \ref{tab.91-3} that $\Phi$ consists of 8 distinct orbits of length 11 and 3 distinct orbits of length 1. Without loss of generality we may assume that the selected orbits of length 1 have 0 in the second component.

We conclude that $s_0 = 11k_0 +3$, $s_1= 11k_1$ and $s_2 =  11k_2$ for some integers $k_0, k_1$ and $k_2$.
Then, by (\ref{eqn.100}) we obtain that $(11k_0+3)^2+(11k_1)^2+(11k_2)^2=2821$, But 2812 is not divisible by 11, which is a contradiction.

\end{pf}

\begin{table}[!h]
\caption{Orbits of $G = \Z_{92} \times \Z_3$ under $x \ra x^{13}$}
\label{tab.91-3}
\begin{center}
\small{
\begin{tabular}{|l|}
\hline
orbits of length 1\\
	  \{(0, 0)\}	\{(0, 1)\}    \{(0, 2)\}\\	
	  \{(23, 0)\}	\{(23, 1)\}    \{(23, 2)\}\\    
	  \{(46, 0)\}	\{(46, 1)\}    \{(46, 2)\}\\
    \{(69, 0)\}    \{(69, 1)\}    \{(69, 2)\}\\
	\hline
	orbits of length 11\\
  \{(33, 2), (45, 2), (61, 2), (5, 2), (17, 2), (57, 2), (89, 2), (65, 2), (37, 2), (21, 2), (53, 2)\}\\
  \{(31, 0), (75, 0), (35, 0), (87, 0), (55, 0), (27, 0), (47, 0), (71, 0), (3, 0), (59, 0), (39, 0)\}\\
\hline
\end{tabular}
}
\end{center}
\end{table}

\section*{Acknowledgment}
The author is supported by the Scientific and Technological Research Council of Turkey (T\"UB\. ITAK) under the National Postdoctoral Research Scholarship NO 2219. 

\bibliographystyle{splncs_srt}
\bibliography{DPDS-gamma}

\begin{thebibliography}{10}

\bibitem{BJL}
Beth, T., Jungnickel, D., Lenz, H.:
\newblock Design theory. {V}ol. {I}. Second edn. Volume~69 of Encyclopedia of
  Mathematics and its Applications.
\newblock Cambridge University Press, Cambridge (1999)

\bibitem{Brock}
Brock, B.W.:
\newblock Hermitian congruence and the existence and completion of generalized
  {H}adamard matrices.
\newblock J. Combin. Theory Ser. A \textbf{49}(2) (1988)  233--261

\bibitem{CTZ2010}
Chee, Y.M., Tan, Y., Zhou, Y.:
\newblock Almost {$p$}-ary perfect sequences.
\newblock In Carlet, C., Pott, A., eds.: Sequences and their
  applications---{SETA} 2010. Volume 6338 of Lecture Notes in Comput. Sci.
\newblock Springer, Berlin (2010)  399--415

\bibitem{HK1998}
Helleseth, T., Kumar, P.V.:
\newblock Sequences with low correlation.
\newblock In Pless, V.S., Huffman, W.C., Brualdi, R.A., eds.: Handbook of
  coding theory, {V}ol. {I}, {II}.
\newblock North-Holland, Amsterdam (1998)  1765--1853

\bibitem{JuPo99}
Jungnickel, D., Pott, A.:
\newblock Perfect and almost perfect sequences.
\newblock Discrete Appl. Math. \textbf{95}(1-3) (1999)  331--359

\bibitem{MaSc1995}
Ma, S.L., Schmidt, B.:
\newblock On {$(p^a,p,p^a,p^{a-1})$}-relative difference sets.
\newblock Des. Codes Cryptogr. \textbf{6}(1) (1995)  57--71

\bibitem{MaNg2009}
Ma, S.L., Ng, W.S.:
\newblock On non-existence of perfect and nearly perfect sequences.
\newblock Int. J. Inf. Coding Theory \textbf{1}(1) (2009)  15--38

\bibitem{MaPo1995}
Ma, S.L., Pott, A.:
\newblock Relative difference sets, planar functions, and generalized
  {H}adamard matrices.
\newblock J. Algebra \textbf{175}(2) (1995)  505--525

\bibitem{OYY2012}
{\"O}zbudak, F., Yayla, O., Y{\i}ld{\i}r{\i}m, C.C.:
\newblock Nonexistence of certain almost {$p$}-ary perfect sequences.
\newblock In Helleseth, T., Jedwab, J., eds.: Sequences and their
  applications---{SETA} 2012. Volume 7280 of Lecture Notes in Comput. Sci.
\newblock Springer, Heidelberg (2012)  13--24

\bibitem{pott1995}
Pott, A.:
\newblock Finite geometry and character theory. Volume 1601 of Lecture Notes in
  Mathematics.
\newblock Springer-Verlag, Berlin (1995)

\bibitem{Tur1965}
Turyn, R.J.:
\newblock Character sums and difference sets.
\newblock Pacific J. Math. \textbf{15} (1965)  319--346

\bibitem{WYZ2014}
Winterhof, A., Yayla, O., Ziegler, V.:
\newblock Non-existence of some nearly perfect sequences, near
  {B}utson-{H}adamard matrices, and near conference matrices.
\newblock arXiv preprint arXiv:1407.6548 (2014)

\end{thebibliography}

\newpage
\appendix
\section{Tables }

\begin{multicols}{2}

\setbox\ltmcbox\vbox{
\makeatletter\col@number\@ne
\tiny
\begin{longtable}{|l|l|l|}
\caption[m-ary]{$p$-ary nearly perfect sequences of period $n$ and type  $\gamma = 0$.}
\label{tab.m_ary_gamma_0}\\

\hline \multicolumn{1}{|c|}{$n$} & \multicolumn{1}{c|}{$p$} & \multicolumn{1}{c|}{Comments} \\ \hline 
\endfirsthead

\multicolumn{3}{c}%
{{\bfseries \tablename\ \thetable{} -- continued from previous page}} \\
\hline \multicolumn{1}{|c|}{$v$} &
\multicolumn{1}{c|}{$m$} &
\multicolumn{1}{c|}{Comments} \\ \hline 
\endhead

\hline \multicolumn{3}{|r|}{{Continued on next page}} \\ \hline
\endfoot

\hline 
\endlastfoot

2
 &2 &exists by \cite[Theorem 2.2.9]{pott1995}\\
\hline
3
 &3 &exists by \cite[Theorem 2.2.9]{pott1995}\\
\hline
4
 &2 &exists by \cite[Theorem 2.3]{MaSc1995}\\
\hline
5
 &5 &exists by \cite[Theorem 2.2.9]{pott1995}\\
\hline
6
 &2 &not exists by \cite[Theorem 4.7]{MaNg2009} with q=3 \\
\hline
 &3 &not exists by \cite[Corollary 4.10]{MaNg2009} with q=3 \\
\hline
7
 &7 &exists by \cite[Theorem 2.2.9]{pott1995}\\
\hline
8
 &2 &not exists by \cite[Theorem 4.6]{MaNg2009} with q=2 \\
\hline

\end{longtable}

\unskip
\unpenalty
\unpenalty
}

\unvbox\ltmcbox

\end{multicols}

\newpage
\begin{multicols}{2}

\setbox\ltmcbox\vbox{
\makeatletter\col@number\@ne
\tiny
\begin{longtable}{|l|l|l|}
\caption[m-ary]{$p$-ary nearly perfect sequences of period $n$ and type  $\gamma = -1$.}
\label{tab.m_ary_gamma_m1}\\

\hline \multicolumn{1}{|c|}{$n$} & \multicolumn{1}{c|}{$p$} & \multicolumn{1}{c|}{Comments} \\ \hline 
\endfirsthead

\multicolumn{3}{c}%
{{\bfseries \tablename\ \thetable{} -- continued from previous page}} \\
\hline \multicolumn{1}{|c|}{$v$} &
\multicolumn{1}{c|}{$m$} &
\multicolumn{1}{c|}{Comments} \\ \hline 
\endhead

\hline \multicolumn{3}{|r|}{{Continued on next page}} \\ \hline
\endfoot

\hline 
\endlastfoot

2
 &3 &exists by \cite{HK1998}\\
\hline
3
 &2 &exists by \cite[Corollary 2.8]{JuPo99} \\
\hline
4
 &5 &exists by \cite{HK1998}\\
\hline
5
 &2 &not exists by Theorem 2 (ii) with q=3 and u=5 \\
\hline
 &3 &exists and given in \cite{MaNg2009} \\
\hline
6
 &7 &exists by \cite{HK1998}\\
\hline
7
 &2 &exists by \cite[Corollary 2.8]{JuPo99} \\
\hline
8
 &3 &exists by \cite{HK1998}\\
\hline
\end{longtable}

\unskip
\unpenalty
\unpenalty
}

\unvbox\ltmcbox

\end{multicols}

\newpage
\begin{multicols}{2}

\setbox\ltmcbox\vbox{
\makeatletter\col@number\@ne
\tiny
\begin{longtable}{|l|l|l|}
\caption[m-ary]{$p$-ary nearly perfect sequences of period $n$ and type  $\gamma = 1$.}
\label{tab.m_ary_gamma_1}\\

\hline \multicolumn{1}{|c|}{$n$} & \multicolumn{1}{c|}{$p$} & \multicolumn{1}{c|}{Comments} \\ \hline 
\endfirsthead

\multicolumn{3}{c}%
{{\bfseries \tablename\ \thetable{} -- continued from previous page}} \\
\hline \multicolumn{1}{|c|}{$v$} &
\multicolumn{1}{c|}{$m$} &
\multicolumn{1}{c|}{Comments} \\ \hline 
\endhead

\hline \multicolumn{3}{|r|}{{Continued on next page}} \\ \hline
\endfoot

\hline 
\endlastfoot

3
 &2 &not exists by Theorem 2 (iv) with q=5 and u=1 \\
\hline
4
 &3 & \\
\hline
5
 &2 &exists by an exhaustive search \\
\hline
6
 &5 &not exists by Theorem 2 (iii) with q=5 and u=6 \\
\hline
7
 &2 &not exists by Theorem 2 (ii) with q=3 and u=7 \\
\hline
 &3 & \\
\hline
8
 &7 &not exists by Theorem 2 (iii) with q=7 and u=8 \\
\hline
9
 &2 &not exists by Theorem 2 (iv) with q=17 and u=1 \\
\hline

\end{longtable}

\unskip
\unpenalty
\unpenalty
}

\unvbox\ltmcbox

\end{multicols}

\newpage
\begin{multicols}{2}

\setbox\ltmcbox\vbox{
\makeatletter\col@number\@ne
\tiny
\begin{longtable}{|l|l|l|}
\caption[m-ary]{$p$-ary nearly perfect sequences of period $n$ and type  $\gamma = 2$.}
\label{tab.m_ary_gamma_2}\\

\hline \multicolumn{1}{|c|}{$n$} & \multicolumn{1}{c|}{$p$} & \multicolumn{1}{c|}{Comments} \\ \hline 
\endfirsthead

\multicolumn{3}{c}%
{{\bfseries \tablename\ \thetable{} -- continued from previous page}} \\
\hline \multicolumn{1}{|c|}{$v$} &
\multicolumn{1}{c|}{$m$} &
\multicolumn{1}{c|}{Comments} \\ \hline 
\endhead

\hline \multicolumn{3}{|r|}{{Continued on next page}} \\ \hline
\endfoot

\hline 
\endlastfoot

4
 &2 &not exists by Theorem 2 (iv) with q=5 and u=1 \\
\hline
5
 &3 &exists and given in Section 3 \\
\hline
6
 &2 &exists an exhaustive search\\
\hline
7
 &5 &not exists by Theorem 2 (iii) with q=5 and u=7 \\
\hline
8
 &2 &not exists by Theorem 2 (ii) with q=3 and u=2 \\
\hline
 &3 &not exists by Theorem 2 (iv) with q=2 and u=1 \\
\hline
9
 &7 &not exists by an exhaustive search  \\
\hline
10
 &2 &not exists by Theorem 2 (iv) with q=7 and u=1 \\
\hline
\end{longtable}

\unskip
\unpenalty
\unpenalty
}

\unvbox\ltmcbox

\end{multicols}

\newpage
\begin{multicols}{2}

\setbox\ltmcbox\vbox{
\makeatletter\col@number\@ne
\tiny
\begin{longtable}{|l|l|l|}
\caption[m-ary]{almost $p$-ary nearly perfect sequences of period $n+1$ and type  $\gamma = 0$.}
\label{tab.almost_m_ary_gamma_0}\\

\hline \multicolumn{1}{|c|}{$n$} & \multicolumn{1}{c|}{$p$} & \multicolumn{1}{c|}{Comments} \\ \hline 
\endfirsthead

\multicolumn{3}{c}%
{{\bfseries \tablename\ \thetable{} -- continued from previous page}} \\
\hline \multicolumn{1}{|c|}{$v$} &
\multicolumn{1}{c|}{$m$} &
\multicolumn{1}{c|}{Comments} \\ \hline 
\endhead

\hline \multicolumn{3}{|r|}{{Continued on next page}} \\ \hline
\endfoot

\hline 
\endlastfoot
3
 &2 &not exists by \cite[Theorem 2]{CTZ2010}\\
\hline
4
 &3 &exists by \cite[Example 1]{CTZ2010}\\
\hline
5
 &2 &not exists by \cite[Theorem 2]{CTZ2010}\\
\hline
6
 &5 &not exists by \cite[Result 5]{CTZ2010}\\
\hline
7
 &2 &not exists by \cite[Theorem 2]{CTZ2010}\\
\hline
 &3 &exists by \cite[Example 1]{CTZ2010}\\
\hline
8
 &7 &exists by \cite[Example 1]{CTZ2010}\\
\hline
9
 &2 &not exists by \cite[Theorem 2]{CTZ2010}\\
\hline
\end{longtable}

\unskip
\unpenalty
\unpenalty
}

\unvbox\ltmcbox

\end{multicols}

\newpage
\begin{multicols}{2}

\setbox\ltmcbox\vbox{
\makeatletter\col@number\@ne
\tiny
\begin{longtable}{|l|l|l|}
\caption[m-ary]{almost $p$-ary nearly perfect sequences of period $n+1$ and type  $\gamma = -1$.}
\label{tab.almost_m_ary_gamma_m1}\\

\hline \multicolumn{1}{|c|}{$n$} & \multicolumn{1}{c|}{$p$} & \multicolumn{1}{c|}{Comments} \\ \hline 
\endfirsthead

\multicolumn{3}{c}%
{{\bfseries \tablename\ \thetable{} -- continued from previous page}} \\
\hline \multicolumn{1}{|c|}{$v$} &
\multicolumn{1}{c|}{$m$} &
\multicolumn{1}{c|}{Comments} \\ \hline 
\endhead

\hline \multicolumn{3}{|r|}{{Continued on next page}} \\ \hline
\endfoot

\hline 
\endlastfoot

2
 &2 &exists by \cite[Example 3]{CTZ2010}\\
\hline
3
 &3 &not exists by \cite[Theorem 7]{CTZ2010}\\
\hline
4
 &2 &exists by \cite[Example 3]{CTZ2010}\\
\hline
5
 &5 &not exists by \cite[Theorem 7]{CTZ2010}\\
\hline
6
 &2 &exists by \cite[Example 3]{CTZ2010}\\
\hline
 &3 &exists by \cite[Example 3]{CTZ2010}\\
\hline
7
 &7 &not exists by \cite[Theorem 7]{CTZ2010}\\
\hline
8
 &2 &not exists by Theorem 4 (ii) with q=3 and u=9 \\
\hline

\end{longtable}

\unskip
\unpenalty
\unpenalty
}

\unvbox\ltmcbox

\end{multicols}

\newpage
\begin{multicols}{2}

\setbox\ltmcbox\vbox{
\makeatletter\col@number\@ne
\tiny
\begin{longtable}{|l|l|l|}
\caption[m-ary]{almost $p$-ary nearly perfect sequences of period $n+1$ and type  $\gamma = 1$.}
\label{tab.almost_m_ary_gamma_1}\\

\hline \multicolumn{1}{|c|}{$n$} & \multicolumn{1}{c|}{$p$} & \multicolumn{1}{c|}{Comments} \\ \hline 
\endfirsthead

\multicolumn{3}{c}%
{{\bfseries \tablename\ \thetable{} -- continued from previous page}} \\
\hline \multicolumn{1}{|c|}{$v$} &
\multicolumn{1}{c|}{$m$} &
\multicolumn{1}{c|}{Comments} \\ \hline 
\endhead

\hline \multicolumn{3}{|r|}{{Continued on next page}} \\ \hline
\endfoot

\hline 
\endlastfoot

4
 &2 &not exists by Theorem 4 (iii) resp. with q=3 and u=5 \\
\hline
5
 &3 &not exists by Theorem 4 (iii) with q=2 and u=1 \\
\hline
6
 &2 &not exists by Theorem 4 (iii) with q=3 and u=1 \\
\hline
7
 &5 &not exists by Theorem 4 (iii) with q=2 and u=1 \\
\hline
8
 &2 &not exists by an exhaustive search \\
\hline
 &3 &not exists by an exhaustive search \\
\hline
9
 &7 &not exists by an exhaustive search \\
\hline
10
 &2 &not exists by Theorem 4 (iii) with q=5 and u=1 \\
\hline

\end{longtable}

\unskip
\unpenalty
\unpenalty
}

\unvbox\ltmcbox

\end{multicols}

\newpage
\begin{multicols}{2}

\setbox\ltmcbox\vbox{
\makeatletter\col@number\@ne
\tiny
\begin{longtable}{|l|l|l|}
\caption[m-ary]{almost $p$-ary nearly perfect sequences of period $n+1$ and type  $\gamma = 2$.}
\label{tab.almost_m_ary_gamma_2}\\

\hline \multicolumn{1}{|c|}{$n$} & \multicolumn{1}{c|}{$p$} & \multicolumn{1}{c|}{Comments} \\ \hline 
\endfirsthead

\multicolumn{3}{c}%
{{\bfseries \tablename\ \thetable{} -- continued from previous page}} \\
\hline \multicolumn{1}{|c|}{$v$} &
\multicolumn{1}{c|}{$m$} &
\multicolumn{1}{c|}{Comments} \\ \hline 
\endhead

\hline \multicolumn{3}{|r|}{{Continued on next page}} \\ \hline
\endfoot

\hline 
\endlastfoot

5
 &2 &not exists by Theorem 4 (iii) with q=3 and u=1 \\
\hline
6
 &3 &not exists by Theorem 4 (iii) with q=2 and u=1 \\
\hline
7
 &2 &not exists by Theorem 4 (iii) with q=3 and u=1 \\
\hline
8
 &5 &not exists by Theorem 4 (iii) with q=2 and u=1 \\
\hline
9
 &2 &not exists by Theorem 4 (iii) with q=3 and u=1 \\
\hline
 &3 &not exists by an exhaustive search\\
\hline
10
 &7 &not exists by Theorem 4 (iii) with q=3 and u=1 \\
\hline
11
 &2 &not exists by Theorem 4 (iii) with q=3 and u=1 \\
\hline
\end{longtable}

\unskip
\unpenalty
\unpenalty
}

\unvbox\ltmcbox

\end{multicols}


\end{document}